\input amstex
\input epsf
\magnification = \magstep1
\overfullrule = 0pt

\def\A{{\bold A}}
\def\B{{\bold B}}
\def\C{{\bold C}}
\def\D{{\bold D}}

\def\G{{\bold G}}
\def\I{{\bold I}}
\def\J{{\bold J}}

\def\Q{{\bold Q}}
\def\R{{\bold R}}
\def\T{{\bold T}}

\def\uF{{\bold u_F}}
\def\uI{{\bold u_I}}

\def\e{{\bold e}}

\def\q{{\bold q}}

\def\s{{\bold s}}
\def\x{{\bold x}}
\def\y{{\bold y}}

\def\u{{\bold u}}
\def\vv{{\bold v}}
\def\w{{\bold w}}
\def\One{{\bold 1}}

\def\Zero{\bold 0}

\loadmsbm
\def \Reals{\Bbb R}

\def\G{\Bbb G}
\def\H{\Bbb H}

\parskip = 4 pt

\def\mat(#1;#2;#3;#4)%
{\left(\matrix#1&#2\\ #3&#4\endmatrix\right)}

\def \ma(#1;#2;#3;#4;#5;#6;#7;#8;#9)%
{\left(\matrix#1&#2&#3\\ #4&#5&#6\\ #7&#8&#9\endmatrix\right)}

\vglue -10 pt
\centerline{\bf Quaternions in Three Dimensions}

\smallskip
\centerline{\bf Bob Palais }

\smallskip
\noindent
{\narrower
{\bf Abstract. }
We construct the quaternion algebra [10] ``geometrically'' 
by a three dimensional analogue of the classic two
dimensional geometric description of the complex field.
The algebraic description of the multiplication
operation in three dimensions involves the
addition of one term.
The construction leads to novel methods
for implementing and interpolating rotations
and understanding their topology.
\par}

\medskip
\noindent
The bilinear map from ordered pairs of vectors in 
$\Reals^3$ to quaternions,  $\T: \Reals^3 \times \Reals^3 \to \H$ defined by:
$$\T(\vv,\w)=[\vv \cdot \w, \vv \times \w] \eqno(1)$$ 
is known to be surjective. 
(For example, see [3], Theorem 1.)  Note that in [3]  quaternions are viewed as 
$2 \times 2$, self-adjoint, complex matrices, 
while `vectors' are identified with elements of the subspace of such matrices having trace zero.
However here  will view vectors $\vv,\w$  as elements of $\Reals^3$ and, for the
moment, we will pretend that the quaternion structure
had not been invented
and simply view $\T$ as a map to $\Q=[q,\q] \in \Reals \times \Reals^3 \cong \Reals^4$.
Although $\T$ is not injective, if we
define an equivalence relation,
$$(\vv,\w) \sim (\vv',\w') \Leftrightarrow
\vv \cdot \w = \vv' \cdot \w'  \hbox{ and } \vv \times \w = \vv' \times \w' , \eqno(2)$$
then we get a bijective map 
$\tilde{\T}\overline{(\vv,\w)}=\T(\vv,\w)$ 
defined on the set  
$\G = (\Reals^3 \times \Reals^3)/\sim$ 
of its equivalence classes,  the inverse
images $\T^{-1}[q,\q]$ of $\T$.
It  then becomes natural to ask: \par

\noindent
{ \narrower ``Can we introduce a `multiplication'  operation for 
ordered pairs of vectors that makes the non-zero classes, $\G^*$, a group,
and  in such a way that $\tilde{\T}$ is a multiplicative homomorphism
with some corresponding operation on $\Reals^4-{\Zero}$? \par}

\noindent
Our goal in  this note is to show that this question {\bf  does} have an 
affirmative answer, and moreover we will define the multiplication operation 
by generalizing a very familiar geometric construction.  When we are done, 
we will see that its expression as an operation on
$\Reals \times \Reals^3$ is $[q',\q'][q,\q] \to [q'', \q'']$, where
$$  \hbox{\hglue 82 pt} q''=q q' - \q\cdot \q' \hbox{\ \ and \ \ } \q'' =q\q'+  q'\q + \q' \times \q.    \hbox{\hglue 82pt} (3)$$
This formula is a modern representation of the quaternion multiplication formula that 
first appears in [13].
It is not hard to check that with $\e_1,\e_2,\e_3$ denoting the standard basis for $\Reals^3$,
if we set $\One=[1,\Zero], i = [0,\e_1], j = [0,\e_2], k = [0,\e_3]$, then $(3)$ implies the
defining relations for quaternions, famously written on a bridge by
Hamilton in 1844 [10], 
$$i^2=j^2=k^2=-\One \hbox{ and }  ijk=-\One  $$ 
[3].  It follows that $ij=k=-ji, ~ jk=i=-kj, ~ ki=j=-ik.$
Thus our construction makes $\tilde{\T}$ an isomorphism 
of $\G^*$ and $\H^*$ and
gives a purely three-dimensional way to interpret (and even `invent') 
quaternions. After developing these ideas, we examine in this light the implications
of viewing known constructions that relate quaternions and rotations.

\noindent
The familiar geometric construction  to which we
refer combines equivalence classes
of ordered pairs of points in Euclidean space to define vector addition, and 
combines equivalence classes of ordered pairs 
of vectors in the plane to define complex multiplication.
(Both of these are treated in detail in [1] where they are referred to as `Chasles' Relation'---though [16]
questions this terminology. Similar constructions are also used to construct the ring
of integers  from the positive integers and the rationals from the integers.) 
What we will show is that a straightforward generalization of those well-known 
geometric constructions does indeed make $\tilde \T$ a multiplicative isomorphism.\par 

\noindent
Using $||\vv||^2 \,||\w||^2 =|\vv \cdot \w|^2+ ||\vv \times \w||^2$, if we
define $||\Q||^2 = q^2+||\q||^2$,  we see that 
the equivalence class $\T^{-1}[q,\q]$ consists of pairs of vectors orthogonal to $\q$ 
that are separated by a fixed oriented angle $\arccos({q \over ||\Q||})$ and the product of whose
lengths is $||\Q||$.
 (If $\q=\Zero$, this says $\vv$ and $\w$ are parallel in an arbitrary 
direction.) An important fact is that any two of these equivalence classes `overlap', 
by which we mean that there exists a unit vector $\u$ for which there are 
representatives of both classes having $\u$ as their first element 
and also a unit vector $\vv$ for which there are representatives 
having $\vv$ as second elements.  This makes it possible, 
in analogy  with vector addition and complex multiplication,  to define 
the geometric composition of an ordered pair of our equivalence classes,
$(\overline{(\vv,\w)}, \overline{(\vv',\w')})$ by finding representatives
such that the second element of the representative of the first class is equal to the 
first element of the representative of the second, and then `merging' them into
$\overline{(\vv'',\w'')} =\overline{(\vv',\w')}\circ \overline{(\vv,\w)}$, the
class of the first element of the first and the second of the second. 
Of course we must show that this geometric merging
operation $\circ$, is well-defined and that it 
makes the non-zero classes a group 
isomorphic to the multiplicative group of non-zero quaternions.  
Both of these facts are a consequence of the following theorem. 

\smallskip
\noindent
{\bf Theorem 1.} {\it  Let $\vv,\w, \vv', \w'$ be vectors in $\Reals^3$ with $\w = \vv'$ and $||\w||=1$.
Then
$$\vv \cdot \w' =(\vv \cdot \w) (\vv' \cdot \w') - (\vv \times \w)\cdot (\vv' \times \w'), \hbox{\hglue 105pt} (4a)$$
and
$$\vv \times \w' =(\vv \cdot \w) (\vv' \times \w')+  
(\vv' \cdot \w') (\vv \times \w) + (\vv' \times \w')\times (\vv \times \w).  \,\,\,\ (4b)$$ \/}

\noindent 
{\bf Remarks: } Before giving the proof, we note that 
if we put $\T(\vv,\w)=[q,\q]$, $\T(\vv',\w')=[q',\q']$, and $\T(\vv,\w')=[q'',\q'']$ 
the conclusion becomes equation $(3)$. Also note that the restriction $||\w||=1$
may be removed by an appropriate scaling that 
puts $(4a,b)$ in a form that highlights their interpretation 
as a novel vector identity involving just three vectors:

\noindent
{\bf Corollary.}
For any three vectors $\A, \B, \C \in \Reals^3$,
$$(\A \cdot \C)(\B \cdot \B) =(\A \cdot \B) (\B \cdot \C) - (\A \times \B)\cdot (\B \times \C), \hbox{\hglue 105pt} (4a')$$
$$(\A \times \C)(\B \cdot \B) =(\A \cdot \B) (\B \times \C)+  
(\B \cdot \C) (\A \times \B) + ((\B \times \C)\times (\A \times \B)).  \,\,\,(4b')$$

\noindent
While there appears to be an asymmetry, with $\B$ plays a distinguished role,
$\A,\B,$ and $\C$ can be cyclically permuted.
We will see that these remarkable vector identities contain the composition laws for (i)the quaternions, $\H$, 
(ii)the unit quaternions, $SU(2)$, (iii)the complex numbers,  $\C$, 
(iv)the plane rotation group $SO(2)$, and (v)the three-dimensional rotation group, $SO(3)$!

\smallskip
\noindent
{\bf Proof:} 
The identities can be validated by direct brute force comparison
of the two sides, but since our purpose is to derive them, we prefer a geometric
approach. 
By the orientation preserving 
Givens' $\Q\R$ factorization algorithm for the matrix with
columns $\{\vv,\w,\w'\}$, we can always construct 
a positively oriented orthonormal ordered 
basis $\{\e_1,\e_2,\e_3\}$ for $\Reals^3$
such that $\e_1=\w$ and $\vv$ is in the $\e_1-\e_2$-plane. 
To motivate $(4a,b)$we first assume $\w'$ is co-planar with $\vv$ and $\w$, i.e., also in the 
$\e_1-\e_2$-plane. In coordinates with respect to this basis, $\q=<0,0,||\q||>,\q'=<0,0,||\q'||>$
so $\vv= < q, -||\q||,0>$, $\w=<1,0,0>$ and $\w'=< q', ||\q'||,0>$.
We confirm by direct calculation that
$$q''=\vv \cdot \w'=q q' - \q\cdot \q', ~~\q'' =\vv \times \w'=q\q'+  q'\q   \hbox{\hglue 90pt} (5)$$
which we recognize from the geometric meaning of dot and cross products as complex multiplication:
$[\cos(\theta)\cos(\theta')-\sin(\theta)\sin(\theta'),<0,0,\cos(\theta)\sin(\theta')+\cos(\theta')\sin(\theta)>].$
We now let 
$\w'$ be arbitrary. Then, since $\q'$ remains orthogonal to $\w=\e_1$,
it takes the form
$\q'=<0,||\q'||\cos(\phi),||\q'|| \sin(\phi)>$
and so $\w' = < q', ||\q'||\sin(\phi),- ||\q'|| \cos(\phi) >$.
Again by direct computation, 
$$q''=\vv \cdot \w'=q q' - \q\cdot \q', ~~\q'' =\vv \times \w'=q\q'+  q'\q  + <||\q|| ~ ||\q'|| \cos(\phi),0,0>.  \hbox{\hglue 10pt}(6)$$
To complete the derivation and proof of $(4a,b)$ 
we observe that the given data must determine this cross-product in terms 
of equivariant quantities and recognize the last term of $(5)$ as the non-commutative
quantity $\q' \times \q$. 
Since in $(4a,b)$ the dot-product is invariant ($\R\x \cdot \R\y=\x \cdot \y$ ) and the cross-product 
equivariant  ($\R\x \times \R\y=\R(\x \times \y)$), they remain valid not only when
the basis or triple of vectors are simultaneously rotated, but also under any rotation
of the pairs $(\vv,\w)$ and $(\vv',\w')$ within their respective equivalence classes. 
This completes the derivation and proof of $(4a,b)$, and guarantees that our
geometric merging operation $\circ$ is well defined on $\G$, and that
with respect to this operation, $\tilde{\T}$ is a multiplicative isomorphism.
{\bf QED}. 

\smallskip
\noindent
What we are doing here is viewing 
quaternion multiplication as the expression of the
relationships $(4a,b)$ that hold for the dot-product and cross-products
of each pair of a triple of vectors in $\Reals^3$.
(Surprisingly, $(4a,b)$ do not appear in standard collections
of similar vector identities.)
While this seems to be novel, 
if we restrict the vectors in $\G$ to lie in any two-dimensional
subspace, the final term in $(4b),(4b')$ vanishes and
$\tilde{\T}$ becomes equivalent to Wessel's
construction of an isomorphism between
the nonzero elements of
$(\Reals^2 \times \Reals^2)/\sim$ and the multiplicative group
of complex numbers, $\C -\{0\}$. 
We obtain a novel description of the three-sphere by
restricting the vectors in $\G$ to be unit vectors, in which case
the $(\B \cdot \B)$ factors in $(4a'),(4b')$ vanish, and
$\tilde{\T}$ becomes an isomorphism between
$(S^2 \times S^2)/\sim$ and the unit quaternions, $SU(2) \cong S^3$.
In this case we may also associate ordered pairs of unit vectors with oriented arcs of
great circles in the unit sphere $S^2 \subset \Reals^3$,
in which case multiplication 
of unit quaternions may be viewed as `vector addition' of these arcs.  
The unit quaternion inverse to $(\vv,\w)$ corresponds to $(\w,\vv)$ and 
the oppositely oriented arc. Algebraically,
the dot-product stays the same while the cross-product is multiplied by $-1$. 
Under the combination of the above restrictions,  the final term in $(4b),(4b')$
and the $(\B \cdot \B)$ factors in $(4a'),(4b')$ vanish, in which case
they revert to the standard trigonometric addition formulas, and
$\tilde{\T}$ becomes an isomorphism between
$(S^1 \times S^1)/\sim$ and the unit circle group, $U(1)\cong SO(2)$. \par 


Continuing the analogy with vectors and complex numbers,
we note that the additive structure
of quaternions is also captured geometrically by $\G$,
using the fact noted earlier that any two equivalence classes
overlap with representatives
of the corresponding equivalence classes that have the same first elements:
$\Q=\T(\vv,\w), \Q'=\T(\vv,\w')$.
Linear combinations of the second element are preserved by $\T$ because of
the bilinearity of both dot-products and cross-products, i.e., 
$c\Q+c'\Q' =\T(\vv,c\w+c'\w')$.
In particular a convex combination (and so  interpolation) 
is preserved.  But even more is true ! We can show that
`SLERP' (quaternion) interpolation [14] on $S^3$: 
$$\Q(t)=c(t)\Q+c'(t)\Q', ~~c(t) = { \sin ((1-t)\Omega)  \over \sin(\Omega) } , ~c'(t)= { \sin(t\Omega) \over \sin(\Omega) }, 
 ~\Omega = \arccos(\Q \cdot \Q') \eqno(7)$$
{\bf can be performed entirely on} $S^2$. To see this, 
note that $\Q(t)=\T(\vv, c(t)\w+c'(t)\w')$
with the same $c(t), ~c'(t)$ as in $(7)$.
Then since 
$$\Q \cdot \Q' =qq'+\q \cdot \q' = (\vv \cdot \w )(\vv \cdot \w' )+ (\vv \times \w )\cdot (\vv \times \w' ),$$
the Lagrange-Binet-Cauchy identity
$$(\A \times \B)\cdot (\C \times \D)= (\A \cdot \C)(\B \cdot \D) - (\A \cdot \D)(\B \cdot \C) \eqno(8)$$
reduces this to $\Q \cdot \Q' = \w \cdot \w'$.  In other words $c(t)$ and $c'(t)$ do not change if we
replace the definition of $\Omega$ in $(7)$ with $\Omega = \arccos(\w \cdot \w')$.
 
 \smallskip
\noindent
This viewpoint also provides insight into the well-known double-covering homomorphism 
of $SU(2)$, the group of unit quaternions, onto the group $SO(3)$ of
rotations of $\Reals^3$.  This takes $\Q \in SU(2)$ to  
$\R(\Q)\in SO(3)$
defined by conjugation:
$[c,\vv] \mapsto \Q [c,\vv]\Q^{-1} = [c,\R(\Q)\vv]$, or equivalently by the
Euler-Rodrigues-Shoemake formula [9]:
$$\R(\Q)= \I  + 2  q \J_\q +  2 \J^2_\q , \hbox{\ where\ } 
 \J_{\u} := \ma ( 0; -u_3; u_2;  u_3 ; 0 ; -u_1 ; -u_2 ; u_1; 0). \eqno(9)$$
When we view a unit quaternion as
$\overline{(\vv,\w)} \in \G$, this homomorphism take the simple form 
$\tilde{C}:\overline{(\vv,\w)} \ \mapsto C(\vv,\w) =
\rho_\w \rho_\vv,$ where $\rho_\u = 2\u \u^T - \I$
is the orthogonal reflection across the one-dimensional subspace
spanned by the unit vector $\u$. The fact that the double covering map $\tilde{C}$
is well-defined on equivalence classes of $(2)$,
is a consequence of the relationship between  rotations $\R$ that fix $\q=\vv\times \w$,  and
reflections $\rho_\u$ with $\u \cdot \q=0$:
$$\rho_{{}_{\R\u}}=\R \,\rho_\u \,\R^{-1} \ .\eqno(10)$$
Then the fact that $\tilde{C}$ is a homomorphism is a consequence of  the
idempotency of the reflections in the unit vector, $\u$, used to compose
two classes, $(\rho_\u)^2 = \I$. 
Because $\rho_{-\u}=\rho_\u$, 
$C(\vv,\w) =C(-\vv,\w)$
where $(-\vv,\w) \sim (\vv,-\w)$ both represent the negative
of $\overline{(\vv,\w)}$. As oriented arcs, they are represented by the
arc complementary to that of $(\vv,\w)$, from which $2(\theta+\pi)=2\theta \mod 2\pi$
also shows that they correspond to the same rotation.
The remaining fact, that its kernel consists {\bf only} 
of the two classes  $(\u,\u)$ and $(\u,-\u)$ for any $\u \in S^2$ follows
from the observation that $\rho_\u \u = \u$ and $\rho_\vv \vv = \vv$
and therefore if $\vv$ is not a multiple of $\u$, 
$C(\u,\vv)$ cannot
fix either $\u$ or $\vv$. 
Therefore we may view $SO(3)$ as equivalence classes of lines
though the origin, with a geometric composition analogous to $\circ$,
and an algebraic composition derived from $(4a),(4b)$ by identifying
classes with opposite signs.

\par 

\noindent
The map $C$ is related to two algorithms for implementing the
unique rotation $\R$ that takes a unit vector $\uI$ to another unit vector $\uF\neq - \uI$
and fixes the orthogonal complement of their span.
The first, motivated by Cartan's `transvection' isometry [2], has been used for an
especially efficient implementation of
a virtual trackball [5] in the mathematical visualization projects, 3D-XplorMath [11] and the
Virtual Math Museum [12]. It defines $\s=\uI+\uF$ and with implicit normalization
to avoid square roots, then applies $\R^{-1}$ to the vectors $\e$ of
the viewing frame in the form $\rho_\s( \rho_\uF \e)$.
This highlights the fact that the angle between any vector and its image 
under the rotation $C(\vv,\w)$ is {\it twice\/} the angle between $\vv$ and $\w$.
This observation
makes it possible to compose two three-dimensional rotations ``with our fingers'',
by composing the corresponding unit quaternions, thus
forming a virtual $SO(4)$ (or, depending on interpretation, $SU(2)$) ``slide rule''.
If you represent the two planes and half-angles of the rotations using the thumbs and
forefingers of your two hands, and then find the equivalent configuration so that the
interior thumb and forefinger overlap, then the exterior thumb and forefinger
represent the plane and half-angle of the composed rotation!
It also explains the difference between the quarter-turn rotation about $\e_3$
corresponding  to the complex $i$ and the half-turn rotation about $\e_1$
corresponding to the quaternion $i$. 
The other algorithm explicitly calculates the matrix for $\R$ in a novel form:
$$\R= \I - {2\over (\s \cdot \s) }\s\s^T + 2\uF  \uI^T .  \eqno{(11)}$$
This can be implemented using only $18$ multiplications and one division, and 
unlike the formula in [7], $(11)$ remains valid in any dimension!

\smallskip
\noindent
Finally, we observe that the intersection of  $S^2$ with 
the family of planes in $\Reals^3$
containing the standard basis vector $\e_1$ whose normals
vary at unit speed along the geodesic from $\e_3$ to $\e_1$
is $$\e(s,t)=\sin(t/4) ( \sin(t/4)\e_1+\cos(t/4)\e_3 )+ \cos(t/4) ( \cos(s)(-\cos(t/4)\e_1+\sin(t/4)\e_3)+\sin(s)\e_2).$$
This explicit  one-parameter family of closed paths on $S^2$ 
(with base-point $\e_1$) 
starts from one turn about the equator
in the $\e_1-\e_2$-plane, namely 
 $\e(s,0)=\cos(s)\e_1+\sin(s)\e_2$, and  shrinks smoothly
 to the constant path $\e(s,2\pi)=\e_1$.
Then $\R(s,t)=C(\e_1, \rho_{\e(s,t)})$
and $\Q(s,t)=\T(\e_1, \rho_{\e(s,t)})$
are explicit homotopies of the paths in  $SO(4)$
and in $S^3$, respectively, between two turns about
a fixed axis ($\e_3$) and the constant path at the identity.
These deformations are the basis of the famous Dirac Belt Trick [15],
Feynman Wine Glass Trick, and the orientation entanglement
relation demonstration described in [6] and visualized in [8] and [12].
All of these express, in different ways, the fact  that
the fundamental group of $SO(4)$ is isomorphic to $Z_2$, but 
the current perspective relates them to
the less mysterious fact that $S^2$ is simply connected,
and also suggests why the same construction
fails for an odd number of turns or if 
the axis is constrained to be fixed during the deformation,
i.e., for plane rotations.


\medskip
\noindent
{\bf References}

\item{[1]} {M. Audin, { \it Geometry\/}, 
Springer, New York, 2003.}

\item{[2]}  {E. Cartan, {\it Le\c cons sur la geometrie des \'espaces de Riemann\/}, Gauthier-Villars, Paris, 1928. }

\item{[3]} {W. Eberlein, The Geometric Theory of Quaternions,
{\it This Monthly.\/} {\bf 9} (1963) 952--954.  

\item{[4]}  {W. R. Hamilton, {On a new species of imaginary quantities connected with a theory of quaternions\/},
{\it Proceedings of the Royal Irish Academy\/}, {\bf 2} (1843) 424--434. }

\item{[5]}  {K. Henriksen, J. Sporring, K. Hornbaek.  Virtual Trackballs Revisited,
{\it  IEEE Transaction on Visualization and Computer Graphics\/}, 
{\bf 2} (2004) 206--216.

\item{[6]}  {C.W. Misner, K.S. Thorne, J.A. Wheeler, {\it Gravitation\/}, 
Freeman, 1973.}

\item{[7]}  {T. M\"oller and  J. F. Hughes, Efficiently Building a Matrix to Rotate One Vector to Another. 
{\it Journal of Graphics Tools\/} {\bf 4} (1999)  1--4.}

\item{[8]} {B. Palais, {Understanding Quaternions Geometrically in $\Reals^3$.\/}, 
{\it LOCI\/}, To appear.}

\noindent
\hskip 5pt
The relevant demos may be found at the following URLs:

\noindent \hskip 5pt
http://www.math.utah.edu/\%7Epalais/quaternion\%5Frepresentatives.html

\noindent \hskip 5pt
http://www.math.utah.edu/\%7Epalais/quaternion\%5Fcomposition.html

\noindent \hskip 5pt
http://www.math.utah.edu/\%7Epalais/quaternion\%5Finterpolation.html

\item{[9]}  {B. Palais and  R. Palais,  Euler's fixed point theorem: The axis of a rotation, {\it J. Fixed Point Theory Appl.\/} {\bf 2} (2007)
215--220. 
}

\item{[10]} R. Palais, Classification of Real Division Algebras, 
{\it This Monthly.\/} {\bf 75} (1964) 366--8.}

\item{[11]} {R. Palais, et.~al., http://3D-XplorMath.org }

\item{[12]} {R. Palais, et.~al., http://VirtualMathMuseum.org }

\item{[13]} {O. Rodrigues,  Des lois g\'eom\'etriques qui r\'egissent les d\'eplacements d'un syst\`eme
solide dans l'\'espace, {\it J. Math\'ematique Pures et Appliqu\'ees.\/} {\bf 5} (1840) 380--440. }

\item{[14]} {K. Shoemake, {Animating Rotation with Quaternion Curves\/}, 
{\it Computer Graphics\/}, {\bf 19}, (1985) 245--254.}

\item{[15]} {V. Stojanoska and O. Stoytchev, Touching the $Z_2$ in Three-Dimensional Rotations,   
{\it Math Magazine 81 \/} {\bf 5} (2008) 345--357. }

\item{[16]} {R. Zenth, {Chasles' Relation.\/}, 
http://iml.univ-mrs.fr/\%7Eritzenth/divers-loisirs/chasles2.pdf.}

\end
\vfill\eject

\smallskip
\noindent
Observe that we have {\it derived } the algebraic quaternion operation $(3)$ 
as the natural extension of the geometric complex multiplication operation
in the plane to {\it three} dimensions, resulting in the addition of one additional
term to its algebraic form $(5)$. 
Visual examples of this construction and others below may be found online at $[6]$.